\documentclass{amsart}
\usepackage{latexsym,amsmath,amstext,mathrsfs,amssymb,amsthm,amscd,amsopn,verbatim,graphics,amsrefs}

\newtheorem{Thm}{Theorem}[section]

\newtheorem{Lem}[Thm]{Lemma}
\newtheorem{Cor}[Thm]{Corollary}

\theoremstyle{remark}
\newtheorem{Rem}[Thm]{Remark}

\theoremstyle{definition}

\newtheorem{Def}[Thm]{Definition}
\newtheorem{Exa}[Thm]{Example}

\newtheorem*{Reg*}{Regularization procedure}

\newcommand{\K}{\mathbb C}
\title[Differential operators on toric varieties]
{Differential operators on toric varieties and Fourier transform}
\author{Giovanni Felder}
\author{Carlo A. Rossi}
\address{Department of mathematics, ETH Zurich, 8092 Zurich, Switzerland}
\email{felder@math.ethz.ch, crossi@math.ethz.ch}

\begin{document}

\begin{abstract}
We show that Fourier transforms on the Weyl algebras have a geometric counterpart in the framework of toric varieties, namely they induce isomorphisms between twisted rings of differential operators on regular toric varieties, whose fans are related to each other by reflections of  one-dimensional cones.
The simplest class of examples is provided by the toric varieties related by such reflections to projective spaces. It includes the blow-up at a point in affine space and resolution of singularities of varieties appearing in the study of the minimal orbit of $\mathfrak{sl}_{n+1}$.
\end{abstract}

\maketitle
\section{Introduction}\label{s-0}
Let $p\colon\widetilde{\mathbb A^n}\to\mathbb A^n=\K^n$ ($n\geq2$)
be the blow-up at the origin of the $n$-dimensional affine space
over the complex numbers (or an algebraically closed field of
characteristic $0$) and let $E=p^{-1}(0)$ be the exceptional
divisor. Consider the ring $\mathcal D_{  \mathcal{O} (m E)}(\widetilde {\mathbb A^n})$ of global
regular differential operators twisted by (i.e., acting on sections
of) $  \mathcal{O} (m E)$, the line bundle associated with the divisor $m E$,
$m\in\mathbb Z$. Our first observation is that this ring is
isomorphic to the ring $\mathcal D_{\mathcal{O}(m-n)}(\mathbb P^n)$ of
differential operators on the $n$-dimensional projective space
twisted by the line bundle $\mathcal{O}(m-n)$. To describe the isomorphism we
realise both rings as subrings of the ring
\[
\mathcal D({\mathbb A^n})=
\K[z_1,\dots,z_n;\partial_1,\dots,\partial_n]
\]
of differential operators in $n$ variables with polynomial
coefficients. Namely any section $\sigma_{mE}$ with divisor $mE$ can
be used to trivialize the bundle over the complement
$\widetilde{\mathbb A^n}\smallsetminus E\simeq\mathbb
A^n\smallsetminus\{0\}$ of $E$ and thus we have an injective
restriction map $ \mathcal D_{  \mathcal{O} (m E)}(\widetilde{\mathbb
A^n})\hookrightarrow \mathcal D(\mathbb A^n\smallsetminus\{0\})$.
Since $\{0\}$ is of codimension at least $2$ in $\mathbb A^n$
(recall that $n\geq 2$), we have an isomorphism $\mathcal D(\mathbb
A^n\smallsetminus\{0\})=\mathcal D(\mathbb A^n)$. Thus we get a ring
monomorphism
\[
i_m\colon \mathcal D_{  \mathcal{O} (m E)}(\widetilde{\mathbb A^n})\hookrightarrow
\mathcal D(\mathbb A^n).
\]
Similarly, the bundle $\mathcal{O}(\ell)$ is associated with, say, a multiple
of the hyperplane in $\mathbb P^n$ defined by the vanishing of some
homogeneous coordinate, whose complement is $\mathbb A^n$. So again
we have an injective restriction map
\[
j_\ell\colon\mathcal D_{\mathcal{O}(\ell)}(\mathbb P^n)
\hookrightarrow \mathcal D({\mathbb A^n}).
\]
Let finally $F\colon \mathcal D(\mathbb A^n)\to \mathcal D(\mathbb
A^n)$, the {\em Fourier transform}, be the ring automorphism such
that $F(z_i)=\partial_i, F(\partial_i)=-z_i$.
\begin{Thm}\label{t-1} Let $m\in\mathbb Z$ and $n\geq 2$. Then the Fourier
transform restricts to a ring isomorphism
\[
 \phi_m\colon\mathcal D_{\mathcal{O}(m-n)}(\mathbb P^n)\to
 \mathcal D_{  \mathcal{O} (m E)}(\widetilde{\mathbb A^n}).
\]
Namely, we have $F\circ j_{m-n}=i_m\circ \phi_m$.
\end{Thm}
Recall that $\mathbb P^n=SL_{n+1}/P$ is a homogeneous space with
parabolic isotropy subgroup $P$ and that $\mathcal{O}(\ell)$ is the
line bundle associated with a character of $P$. Thus the results of
Borho and Brylinski \cite{BorhoBrylinski}*{Theorem 3.8 and Remark
3.9} apply: there is a surjective algebra homomorphism
\[
U(\mathfrak{sl}_{n+1})\to \mathcal D_{\mathcal{O}(\ell)}(\mathbb P^n),
\]
given by the infinitesimal action of the group $SL_{n+1}$ and the
kernel is the annihilator $J_\ell$ of a generalized Verma module.

\begin{Cor}\label{c-1}
There is a surjective algebra homomorphism
\[
U(\mathfrak{sl}_{n+1})\to \mathcal D_{  \mathcal{O} (m E)}(\widetilde{\mathbb A^n})
\]
with kernel $J_{m-n}$.
\end{Cor}
 The $\K$-algebra $\mathcal D_{\mathcal{O}(\ell)}(\mathbb P^n)$, viewed as a
subalgebra of $\mathcal D(\mathbb A^n)$, is generated by $1$ and the
first order operators $\partial_i,
z_i\partial_j+\delta_{ij}(e-\ell), z_i(e-\ell)$ given by the action
of the fundamental vector fields. Here $e=\sum_iz_i\partial_i$
denotes the Euler vector field.

\begin{Cor}\label{c-2}
The $\K$-algebra $\mathcal D_{  \mathcal{O} (m E)}(\widetilde{\mathbb A^n})$, viewed as a
subalgebra of $\mathcal D(\mathbb A^n)$, is generated by $1$ and the
second order operators $z_i,z_i\partial_j,\partial_i(e+m)$,
$i,j=1,\dots,n$.
\end{Cor}

 The algebra of differential operators twisted by a line bundle
acts on the space of global section of the line bundle. By Corollary
\ref{c-1} the space of global sections $\Gamma(\widetilde{\mathbb
A^n},   \mathcal{O} (m E))$ is a module over $U(\mathfrak{sl}_{n+1})$. Let
$\mathfrak{sl}_{n+1}=\mathfrak n_-\oplus \mathfrak h\oplus \mathfrak
n_+$ be the Cartan decomposition into lower triangular, traceless
diagonal and upper triangular matrices. Denote by $M(\lambda)$ the
Verma module with highest weight $\lambda\in\mathfrak h^*$ and by
$L(\lambda)$ the irreducible quotient of $M(\lambda)$, see
\cite{Jantzen}. The fundamental weights are $\varpi_i\colon x\mapsto
x_1+\dots+x_i$, $x=\mathrm{diag}(x_1,\dots,x_{n+1})\in\mathfrak h$.
\begin{Thm}\label{t-2}
The $U(\mathfrak{sl}_{n+1})$-module $\Gamma(\widetilde{\mathbb A^n},
\mathcal{O} (m E))$ is isomorphic to $L(-(m+1)\varpi_n)$ if $m\geq0$
and to $L((m-1)\varpi_n-m\varpi_{n-1})$ if $m<0$.
\end{Thm}

We deduce Theorem \ref{t-1} from a more general theorem (Theorem
\ref{t-3}) on toric varieties, stating that suitably twisted
algebras of differential operators on nonsingular toric varieties
whose fans are related by reflections of one-dimensional cones, see
Fig.~\ref{f-1} and Section \ref{s-2}, are isomorphic via a Fourier
transform. 
The proof of this theorem uses Musson's description \cite{Musson} of
the algebra of differential operators on a toric variety and 
partial Fourier transforms on Weyl algebras as in the work of
Musson and Rueda \cite{Musson-Rueda}.

In this way we obtain several families of toric varieties
with line bundles whose twisted algebras of differential operators
are isomorphic. Even in the simplest case of projective spaces,
that we treat in detail in Section \ref{s-3}, the situation is
rather rich. Additionally to the projective $n$-space and the blow-up at
the origin of the affine $n$-space there are many other varieties related by such
reflection. They include resolutions of singularities of the
varieties considered by Levasseur, Smith and Stafford
\cite{Levasseuretal} in the case of $\mathfrak{sl}_{n+1}$ and Musson \cite{Musson2}. 
Our
results extends and unify their results. For these varieties we also
compute the action of $\mathfrak{sl}_{n+1}$ on the cohomology 
with coefficients in any line bundle and find that cohomology groups
form irreducible modules.
They are described in Theorem \ref{t-5} and Theorem \ref{thm-irr}.
Theorem \ref{t-2} is a
special case of this more general result. 

Let us note that the
problem of studying isomorphisms of algebras of differential
operators is also actively studied in different contexts, mostly for
singular affine varieties, see \cite{BerestWilson} for a recent
review.

 {\bf Acknowledgements.} We thank
Y. Berest for discussions and explanations and I. Musson for 
pointing out that \cite{Musson-Rueda} already contains part
of our results. 
 This work has been partially
supported by the European Union through the FP6 Marie Curie RTN ENIGMA
(Contract number MRTN-CT-2004-5652), by the Swiss National Science
Foundation (grant 200020-105450) and the MISGAM programme of the European
Science Foundation.

\section{Differential operators on toric varieties and Fourier transform}\label{s-1}
\subsection{Line bundles on toric varieties}
We refer to the book \cite{Fulton} for an introduction to toric varieties
and line bundles on them, and use mostly
the same notations. Let $N\simeq \mathbb Z^n$ be a lattice of rank $n$ and
$M=\mathrm{Hom}(N,\mathbb Z)$ the dual lattice. Recall that a fan in
$N$ is a finite collection of strongly convex rational polyhedral
cones $\sigma\subset N_{\mathbb R}=N\otimes_\mathbb Z \mathbb R$
such that each face of a cone in $\Delta$ is a cone in $\Delta$ and
the intersection of cones in $\Delta$ is a face of each. To each fan
$\Delta$ there corresponds a toric variety $X=X(\Delta)$, a normal
algebraic variety with an action of the torus $T_N=N\times_\mathbb Z
\K^\times\simeq (\K^\times)^n$ with a dense orbit $U_0$. The
functions on $U_0$ are spanned by Laurent monomials $\chi_\mu$ with
exponent $\mu\in M$. To each $\sigma\in\Delta$ there corresponds an
invariant affine open set $U_\sigma$ and an orbit closure
$V(\sigma)$ of codimension equal to the dimension of $\sigma$. The orbit
closures
$D_1,\dots, D_d$ associated with the
one-dimensional cones are $T_N$-invariant Weil divisors
and there is an exact sequence of groups (\cite{Fulton}, Proposition on page 63)
\[
 M\to \oplus_{i=1}^d\mathbb Z\, D_i\to A_{n-1}(X)\to 0,
\]
where $A_{n-1}(X)$ is the group of Weil divisors modulo divisors of
rational functions. Let $k_1,\dots, k_d$ be the nonzero lattice
vectors closest to the origin on the one-dimensional cones. The left
map sends $\mu\in M$ to $\sum\langle\mu,k_i\rangle D_i$, which is the divisor
of the function $\chi_\mu$. The sequence is also exact on the
left if the vectors in $\Delta$ span $N_\mathbb R$.

 We will make the following simplifying
assumptions about $\Delta$.

\begin{Def}
A fan $\Delta$ is {\em regular} if
\begin{enumerate}
\item[(i)] Each cone in $\Delta$ is generated by part of a basis of
$N$.
\item[(ii)] The one-dimensional cones in $\Delta$ are generated by
vectors $k_1,\dots,k_d$ spanning the lattice $N$ over $\mathbb Z$.
\end{enumerate}
\end{Def}
Assumption (i) means that we consider non-singular toric varieties
so that the map $\mathrm{Pic}(X)\to\mathrm A_{n-1}(X)$ sending a
line bundle to the class of the divisor of a non-zero rational
section is an isomorphism. Assumption (ii) implies in particular that
$k_i$ span $N$ and thus there is an exact sequence of (free Abelian) groups
\begin{equation}\label{e-xs1}
0\to M\to \sum_{i=1}^d \mathbb Z D_i\to
A_{n-1}(X)\to 0,\qquad A_{n-1}(X)\simeq \mathrm{Pic}(X).
\end{equation}
Under these assumptions, each cone of $\Delta$ is generated by a
subset of the lattice vectors $k_1,\dots,k_d$, called the {\em
generating vectors} of the fan.

\begin{Exa}\label{ex-Piscine Molitor} Let $\epsilon_1,\dots,\epsilon_n$ be a basis
of $N$ and let $\epsilon=\epsilon_1+\dots+\epsilon_n$. Then the fan
in $N$ whose cones are generated by all proper subsets of
$\{\epsilon_1,\dots,\epsilon_n,-\epsilon\}$ gives the toric variety
$\mathbb P^n$ and the fan whose cones are generated by all proper
subsets of $\{\epsilon_1,\dots,\epsilon_n,\epsilon\}$ except
$\{\epsilon_1,\dots,\epsilon_n\}$ gives the toric variety
$\widetilde{\mathbb A^n}$. In Fig.~\ref{f-1} the case $n=2$ is
depicted. \end{Exa}
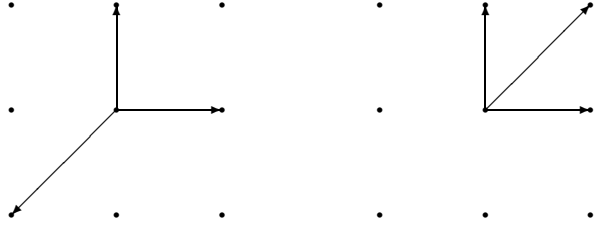
\begin{figure}
\setlength\unitlength{0.7mm}
\begin{picture}(110,40)(0,0)
 \put( 20, 20){\circle*{1}}
 \put( 40, 20){\circle*{1}}
 \put(  0, 20){\circle*{1}}
 \put( 20,  0){\circle*{1}}
 \put( 40,  0){\circle*{1}}
 \put(  0,  0){\circle*{1}}
 \put( 20, 40){\circle*{1}}
 \put( 40, 40){\circle*{1}}
 \put(  0, 40){\circle*{1}}
 \put( 90, 20){\circle*{1}}
 \put(110, 20){\circle*{1}}
 \put( 70, 20){\circle*{1}}
 \put( 90,  0){\circle*{1}}
 \put(110,  0){\circle*{1}}
 \put( 70,  0){\circle*{1}}
 \put( 90, 40){\circle*{1}}
 \put(110, 40){\circle*{1}}
 \put( 70, 40){\circle*{1}}
 \put( 20, 20){\vector(-1,-1){20}}
 \put( 20, 20){\vector( 1, 0){20}}
 \put( 20, 20){\vector( 0, 1){20}}
 \put( 90, 20){\vector( 1, 1){20}}
 \put( 90, 20){\vector( 1, 0){20}}
 \put( 90, 20){\vector( 0, 1){20}}
\end{picture}
\caption{The generating vectors of the fans of $\mathbb P^2$ (left)
and $\widetilde{\mathbb A^2}$ (right)}\label{f-1}
\end{figure}
\subsection{Musson's construction}
Musson \cite{Musson} gave a description of the ring of twisted
differential operators on an arbitrary toric variety in terms of
generators and relations. We recall here his construction in the
special case of varieties with regular fan, for which some
simplifications occur. Let us above $k_1,\dots,k_d\in N$ be the
generating vectors of a regular fan $\Delta$. Let $p\colon \mathbb
Z^d\to N$ be the group homomorphism sending the $i$-th standard
basis vector $\epsilon_i$ to $k_i$, $i=1,\dots,d$. For each
$\sigma\in\Delta$ let $\hat\sigma\in \mathbb R^d$ be the cone
spanned by the vectors $\epsilon_i$ such that $k_i\in\sigma$. The
cones $\hat\sigma$ form a fan in $\mathbb R^d$ and define a toric
variety $Y$ with the action of $(\K^\times)^d$. We have an exact
sequence
\begin{equation}\label{e-xs2}
 0\to K\to \mathbb Z^d\to N\to 0,
\end{equation}
The torus $(\K^\times)^d$ acting on $Y$ has thus a subgroup $G=T_K$ and it is
proved in \cite{Musson} (and also elsewhere, see \cite{Audin}, \cite{Cox}
and references therein)
that $X=Y//G$ is
the algebro-geometric quotient of $Y$ by $G$. The variety $Y$ is the
union of $G$-invariant open sets $V_\sigma$ in $\K^m$ associated
with the cones $\hat\sigma$, $\sigma\in\Delta$:
\[
V_\sigma=\{Q\in \K^d\,|\,Q_i\neq 0, \text{whenever}\ k_i\notin
\sigma\}.
\]
Thus $Y$ is the complement in $\mathbb A^d$ of the union of subspaces of
codimension $\geq 2$.
 Moreover, for each character $\chi\in\mathrm{Hom}(G, \K^\times)$ there
corresponds a sheaf $\mathcal L_\chi$ on $X$: the space of sections
over $U_\sigma=V_\sigma/G$ is the $\mathcal O_X(U_\sigma)=\mathcal
O_Y(V_\sigma)^G$-module
\begin{equation}\label{e-Richard Parker}
\mathcal L_\chi(U_\sigma)=\{f\in O_Y(V_\sigma)\,|\,
f(gx)=\chi(g)f(x), g\in G\}.
\end{equation}
Under our assumptions, $\mathcal L_\chi$ is the sheaf of sections of a line
bundle and each line bundle can be obtained this way. To see this
notice that the exact sequence \eqref{e-xs2} is dual to the exact
sequence \eqref{e-xs1} so
\[\mathrm{Hom}(G,\K^\times)\simeq
\mathrm{Hom}(K,\mathbb Z)\simeq A_{n-1}(X).
\]
Explicitly, the character $\chi(z)=z^a=z_1^{a_1}\cdots z_d^{a_d}$,
$z\in T_K\subset (\K^\times)^d$ corresponds to the class in
$A_{n-1}(X)$ of the divisor $\sum a_iD_i$. The fundamental vector
fields of the infinitesimal action of $G$ on $Y$ form a Lie algebra
$\mathfrak g\simeq K\otimes_{\mathbb Z}\K$ and $\chi$ defines a Lie
algebra homomorphism $\chi_*$, from $\mathfrak g$ to $\K$. The group
$G$ acts on $\mathcal O_Y$ and on $\mathcal D(Y)$. For each
$\sigma\in \Delta$, the restriction of a $G$-invariant differential
operator on $Y$ to $V_\sigma$ maps $\mathcal L_\chi(U_\sigma)$ to
itself, compatibly with inclusions of open sets. Thus there is a map
$\mathcal D(Y)^G\to \mathcal D_{\mathcal L_\chi}(X)$.

\begin{Thm} (Musson, \cite{Musson}*{Theorem 5}) The natural map $\mathcal D(Y)^G\to
\mathcal D_{\mathcal L_\chi}(X)$ is surjective with kernel
generated by $\{x-\chi_*(x)1\,|\, x\in \mathfrak g\}$.
\end{Thm}
Using this theorem, Musson gives a description of $\mathcal
D_{\mathcal L_\chi}(X)$ in terms of generators and relations. Let us
describe it under our assumptions. The complement of $Y$ in $\mathbb
A^d$ has codimension at least $2$ so
\[
\mathcal D(Y)^G\simeq \mathcal D(\mathbb A^d)^G.
\]
The Weyl algebra $\mathcal D(\mathbb A^d)$ is generated by
coordinate functions $Q_i$ and vector fields $P_i=\partial/\partial
Q_i$. A basis is formed by monomials $Q^\lambda
P^\mu=Q_1^{\lambda_1}\cdots P_d^{\mu_d}$, $\lambda,\mu\in\mathbb
N^d$ and $g\in G\subset (\K^\times)^d$ acts on these monomials by
multiplying them by $g^{\lambda-\mu}$. Thus $\mathcal D_{\mathcal
L_\chi}(X)$ is isomorphic to the quotient of the subalgebra of
$\mathcal D(\mathbb A^d)$ generated by $Q^\lambda P^\mu$,
$\lambda-\mu\in K^\perp$ by the ideal generated by $\sum_i
m_iQ_iP_i-\chi_*(m)1$, $m\in K$. For the line bundle
$\mathcal{O}(D)$ associated with the $T_N$-invariant divisor
$D=\sum_{i}a_iD_i$, we have $\chi_*(m)=\sum_{i}a_im_i$.

\begin{Exa}\label{e-Pn} Let $X=\mathbb P^n$ with fan as in Example
\ref{ex-Piscine Molitor}. Then $B=\mathbb Z^{n+1}$ and
$p(\epsilon_i)=\epsilon_i$, $i\leq n$, $p(\epsilon_{n+1})=-e$. The
kernel $K$ is spanned by $(1,\dots,1)$,
$Y=\K^{n+1}\smallsetminus\{0\}$ with coordinates $Q_1,\dots,
Q_{n+1}$, $G=\K^\times$ acts by $Q_i\mapsto \lambda Q_i$,
$\lambda\in\K^\times$ and $\mathfrak g$ is spanned by the vector
field $\sum_{i=1}^{n+1} Q_iP_i$. Then $\mathcal
D_{\mathcal{O}(\ell)}(\mathbb P^n)$ is the quotient of the algebra
of operators $Q^\lambda P^\mu$ ,$\sum_i\lambda_i=\sum_i\mu_i$,
$\lambda_i,\mu_i\geq0$, by the ideal generated by $\sum
Q_iP_i-\ell$.
\end{Exa}

\section{Fourier transforms}\label{s-2}
Let $d\geq n\geq2$ and $I\subset\{1,\dots,d\}$. Let us say that two
regular fans $\Delta$, $\Delta'$ with the same underlying lattice
$N\simeq \mathbb Z^n$ and the same number $d$ of rays ($=$
one-dimensional cones) are {\em related by an $I$-reflection} if
there is a numbering $k_1,\dots,k_d$, $k'_1,\dots,k'_d$ of the
generating vectors of $\Delta$, $\Delta'$ such that
\[
k'_i=
\begin{cases}
-k_i, \text{ if } i \in I, \\
 k_i, \text{otherwise. }
\end{cases}
\]
Let $X=X(\Delta), X'=X(\Delta')$ be the corresponding toric
varieties. The groups $A_{n-1}(X)$, $A_{n-1}(X')$ of classes of Weil
divisors are generated by the classes of the $T_N$-invariant divisors
$D_i$, $D_i'$ associated with the rays $\mathbb R_{\geq0}k_i$,
$\mathbb R_{\geq 0}k_i'$. Let us introduce an {\em affine}
isomorphism $\phi_I\colon A_{n-1}(X)\to A_{n-1}(X')$ defined on
representatives by
\begin{equation}\label{e-Okimoto}
\phi_I\colon\sum_{i=1}^da_iD_i\mapsto\sum_{i\notin
I}a_iD_i'-\sum_{i\in I}(a_i+1)D_i'.
\end{equation}
\begin{Lem} The map $\phi_I$ is well defined, i.e., independent of
the choice of representatives.
\end{Lem}
\begin{proof} The group $A_{n-1}(X)$ is the quotient of
$\sum_{i=1}^d\mathbb Z D_i$ by the image of $M$, embedded via
$\mu\mapsto \sum_{i=1}^d\langle\mu,k_i\rangle D_i$ and similarly for
$A_{n-1}(X')$. Changing representative thus means replacing $a_i$ by
$a_i+\langle\mu,k_i\rangle$. The image under $\phi_I$ changes then
by $\sum_{i\notin I}\langle\mu,k_i\rangle D_i'-\sum_{i\in
I}\langle\mu,k_i\rangle D_i'=\sum_i \langle\mu,k_i'\rangle D_i'$
which belongs to the image of $M$ in $\sum \mathbb Z D_i'$.
\end{proof}

Let $F_I\in\mathrm{Aut}(\mathcal{D}(\mathbb A^d))$ be the
automorphism acting on generators as
\[
F_I(Q_i)=
\begin{cases}
P_i, \text{ if } i \in I, \\
Q_i, \text{ otherwise,}
\end{cases}
\qquad F_I(P_i)=
\begin{cases}
-Q_i, \text{ if } i \in I, \\
P_i, \text{ otherwise. }
\end{cases}
\]
These automorphisms were first considered in this setting 
by Musson and Rueda in \cite{Musson-Rueda}. The next Theorem
is an extension of a result of these authors (\cite{Musson-Rueda}*{Lemma 5.2}).

\begin{Thm}\label{t-3}
Suppose that $\Delta$, $\Delta'$ are regular fans related by an
$I$-reflection for some $I\subset\{1,\dots,d\}$. Let $[D]\in A_{n-1}(X)$
be a Weil divisor and let $D=\sum a_i D_i$ be a $T_N$-invariant representative.
Then the Fourier
transform $F_I$ restricts to an isomorphism $\mathcal D(Y)^{G}\to
\mathcal D(Y')^{G'}$ of the corresponding Musson algebras, which in
turns descends to an algebra isomorphism
\[
\mathcal D_{  \mathcal{O} (D)}(X)\to \mathcal D_{  \mathcal{O} (\phi_I(D))}(X').
\]
of the algebras of differential operators.
\end{Thm}

\begin{proof} 
The first part of the proof is a matter of going through Musson's construction and is taken from \cite{Musson-Rueda}*{Lemma 5.2}. 
The involution $\sigma_I\colon\mathbb Z^d\to\mathbb Z^d$, changing sign
to the coordinates labeled by $I$, maps $K$ to the kernel $K'$ of
the map $p'\colon\mathbb Z^d\to N$ and thus induces an automorphism
$\sigma_I$ of $T_d=(\K^\times)^d$ mapping the subtorus $G$ to $G'$.
This automorphism maps $z$ to $z'$ with $z'_i=z_i^{-1}$ if $i\in I$
and $z'_i=z_i$ otherwise. The torus $T_d$ acts on generators $Q_i$,
$P_i=\partial/\partial Q_i$ of $\mathcal D(\mathbb A^d)$ via $z\cdot
Q_i=z_iQ_i$, $z\cdot P_i=z_i^{-1}P_i$. From this and the formula for
$F_I$ follows that $F_I(z\cdot x)=\sigma_I(z)\cdot F_I(x)$, $z\in T_d$,
$x\in\mathcal D(\mathbb A^d)$. Therefore $F_I$ induces an
isomorphism of the algebras of invariants.

There remains to show that the kernel of the map $\mathcal D(Y)^G\to
\mathcal D_{  \mathcal{O} (D)}(X)$ is mapped to the kernel of the
corresponding map for $X'$. By Musson's theorem this kernel is
generated by the operators $\xi_m=\sum_im_iQ_iP_i-\chi_*(m)1$, where $m$
runs over $K\subset \mathbb Z^d$ and if $D=\sum_i a_i D_i$, then
$\chi_*(m)=\sum_i a_im_i$. Since $m\in K$ if and only if $m'=\sigma_I(m)\in K'$ we
have:
\begin{eqnarray*}
F_I(\xi_m)&=&\sum_{i=1}^dm_i(F_I(Q_iP_i)-a_i)\\
&=&\sum_{i\notin I}m_iQ_iP_i-\sum_{i\in
I}m_iP_iQ_i-\sum_{i=1}^da_im_i\\
&=&\sum_{i=1}^dm_i'Q_iP_i-\sum_{i\notin I}a_im_i-\sum_{i\in
I}(a_i+1)m_i\\
&=&\sum_{i=1}^dm_i'Q_iP_i-\sum_{i\notin I}a_im'_i+\sum_{i\in
I}(a_i+1)m'_i.
\end{eqnarray*}
Thus $F_I$ maps the generators $\xi_m$ of the kernel to the generators $\xi_{m'}$ of
of the kernel of $\mathcal D(Y')^{G'}\to\mathcal D_{  \mathcal{O} (\phi_I(D))}(X')$.
\end{proof}

\begin{Rem}
More generally, this theorem has a version valid for singular toric varieties.
In this case, however, not all Weil divisors are Cartier divisors and thus
they do not all correspond to line bundles.
It can then happen that $\phi_I$ maps a Cartier divisor
to a divisor which is not Cartier.
Suppose first that only assumption (ii) holds, and that both $D$ and $\phi_I(D)$
are Cartier divisors. Then the statement of Theorem \ref{t-3} holds.
If assumption (ii) is not satisfied, then there are factors of $\K^\times$ in
$Y$ and the algebra $\mathcal D(Y)$ has generators
$Q_1,\dots,Q_k,Q_{k+1}^{\pm1},\dots,Q_d^{\pm1}$ and
$P_i=\partial/\partial Q_i$. The Fourier transform $F_I$ is defined
for $I\subset\{1,\dots,k\}$. Theorem \ref{t-3} holds for these subsets $I$
and for Cartier divisors $D$, $\phi_I(D)$.
 \end{Rem}

\bigskip

\section{The case of the projective space}\label{s-3}
\subsection{Varieties related by reflection to the projective
space}\label{ss-31} The projective $n$-dimensional space $\mathbb
P^n=(\mathbb A^{n+1}\smallsetminus \{0\})/\K^\times$ is a toric
variety with generating vectors $k_1,\dots,k_{n+1}\in N\simeq\mathbb
Z^{n}$ such that $\sum_{i=1}^{n+1}k_i=0$. Varieties related to
$\mathbb P^n$ by an $I$-reflection have thus generating vectors
generating $N$ and obeying
\[
\sum_{i\in I}k_i=\sum_{i\notin I}k_i.
\]
The set $S_I$ of regular fans with these generating vectors is a
finite set with a partial order by inclusion. Here are some
interesting examples.
\begin{enumerate}
\item[(a)] {\em Minimal fans.}
For each $I$, $S_I$ has a minimal element: the fan consisting of
just the one-dimensional cones spanned by the vectors $k_i$. The
variety corresponding to this minimal fan is contained as an open
set, with complement of codimension at least 2, in all varieties
$X(\Delta)$, $\Delta\in S_I$.
\item[(b)] {\em Projective spaces.}
The projective space $\mathbb P^n$ has a maximal fan in $S_I$, with
$I=\varnothing$ or $I=\{1,\dots,n+1\}$.
\item[(c)] {\em Blow-ups.}
The blow-up at the origin of $\mathbb A^n$ is obtained from a
maximal fan in $S_I$, for $I$ a set with $r=1$ or $r=n$ elements.
\item[(d)]
{\em Matrices of rank one and of rank $\leq 1$.} Let
$I=\{1,\dots,r\}$, $2\leq r\leq n-1$. A variety with regular fan in
$S_I$ is $Z_r=(\K^r\smallsetminus
\{0\})\times(\K^{n+1-r}\smallsetminus \{0\})/\K^\times$, where
$z\in\K^\times$ acts by $z\cdot(x,y)=(zx,z^{-1}y)$. It may be
identified via $(x,y)\mapsto x^Ty$ with the variety of $r$ by
$n+1-r$ matrices of rank 1, arising in the study of the orbit
$O_{\mathrm{min}}$ of the highest weight vector in the adjoint
representation of $\mathfrak{sl}_{n+1}$, studied in
\cite{Levasseuretal}. The corresponding fan consists of the cones
generated by $\{k_i\,|\,i\in J\}$ for a subset
$J\subset\{1,\dots,n+1\}$ whose complement contains at least an
element in $I$ and at least an element not in $I$.  The closure
$\bar Z_r=Z_r\cup\{0\}$ in the space of $r$ by $n+1-r$ matrices is a
singular affine toric variety. Its fan (it is not in $S_I$) is
obtained from the fan of $Z_r$ simply by adding the cone generated
by $k_1,\dots,k_{n+1}$. These varieties are the irreducible
components of $\bar O_{\mathrm{min}}\cap \mathfrak{n}_+$, see
\cite{Levasseuretal}.
\item[(e)] {\em Resolution of singularities of the above.}
Let again $I=\{1,\dots,r\}$, $2\leq r\leq n-1$. Let $C_i$ be the
$n$-dimensional cone spanned by the basis $(k_j)_{j\neq i}$. Let
$\Delta^+_I$, resp.\ $\Delta^-_{I}$, be the fan consisting of the cones $C_i$ with $i\in I$, resp.\ $i\notin I$,
and their corresponding faces. One can check that both are regular maximal fans.
Both are subdivisions of the fan of $Z_r$ and thus, by the results of Section 2.6.\ of~\cite{Fulton}, the
corresponding toric varieties $\tilde Z^+_{r}$ and $\tilde Z^-_{r}$ give resolutions of
singularity $\tilde Z^\pm_{r}\to \bar Z_r$.
\end{enumerate}

For any $X=X(\Delta)$, $\Delta\in S_I$, the group $A_{n-1}(X)\simeq
\mathrm{Pic}(X)$ of Weil divisors modulo linear equivalence is the
quotient of the group $\sum \mathbb Z D_i$ of $T_N$-invariant Weil
divisor by the relations $D_i= D_j$ if either $i,j\in I$ or
$i,j\notin I$, and $D_i=-D_j$ otherwise. Thus $A_{n-1}(X)\simeq
\mathbb Z$ generated by any $D_i$.

\begin{Thm}\label{t-4}
Let $X=X(\Delta)$, $\Delta\in S_I$ be an $n$-dimensional toric
variety whose fan is regular with $n+1$ generating vectors obeying
$\sum_{i\in I}k_i=\sum_{i\neq I}k_i$, for some subset
$I\subset\{1,\dots,n+1\}$ with $|I|$ elements. Let $D_0$ be a
divisor linearly equivalent to $D_i$, $i\notin I$ or equivalently to
$-D_i$, $i\in I$. Then for all $\ell\in\mathbb Z$, $F_I$ induces an
isomorphism
\[
 \mathcal D_{\mathcal{O}(\ell)}(\mathbb P^n)
 \longrightarrow
 \mathcal D_{  \mathcal{O} ((\ell+|I|)D_0)}(X).
\]
\end{Thm}

\begin{proof}
This is the special case of Theorem \ref{t-3} in which one of the
varieties is $\mathbb P^n$. The line bundle $\mathcal{O}(\ell)$ is isomorphic
to $  \mathcal{O} (\ell D_i)$ for any $i$. Take $i\in I$ (if $I$ is empty there is
nothing to prove). The map $\phi_I$ of Theorem \ref{t-3} sends $\ell
D_i$ to $-\ell D_i-\sum_{j\in I} D_j$, see \eqref{e-Okimoto} which is
equivalent to $-(\ell+|I|)D_i$.
\end{proof}

If $|I|=\{1,\dots,n\}$, $X$ is the blow-up $\widetilde{\mathbb A^n}$, the
exceptional divisor $E$ is $D_{n+1}$ and we obtain Theorem
\ref{t-1}. The formula in affine coordinates for $F_I$ is then the
isomorphism $F$ of the Introduction in this case.

\subsection{The module of global sections}
Let $X$, $\Delta$, $I$ and $D_0$ as in Theorem \ref{t-4}. The
algebra $\mathcal D_{  \mathcal{O} (m D_0)}(X)$, $\Delta\in S_I$,
acts on the space $\Gamma(X,  \mathcal{O} (m D_0))$ of global
sections. Let $\ell=m-|I|$. The homomorphism
$U(\mathfrak{sl}_{n+1})\to \mathcal D_{\mathcal O(\ell D_0)}(\mathbb
P^n)$ composed with the isomorphism of Theorem \ref{t-4} gives
$\Gamma(X, \mathcal{O} (m D_0))$ the structure of an
$\mathfrak{sl}_{n+1}$-module. Let $\varpi_i\colon x\mapsto
x_1+\dots+x_i$, $x=\mathrm{diag}(x_1,\dots,x_{n+1})\in\mathfrak h$
be the fundamental weights of $\mathfrak{sl}_{n+1}=\mathfrak
n_-\oplus\mathfrak h\oplus \mathfrak n_+$ with respect to the
standard Cartan decomposition. Let $s_i\in\mathrm{End}(\mathfrak
h^*)$, $i=1,\dots,n$, denote the simple reflections in the Weyl
group.

\begin{Thm}\label{t-5}
Let $I=\{1,\dots,r\}$ and $X$, $D_0$ be as in Theorem \ref{t-4}. The
$\mathfrak{sl}_{n+1}$-module $\Gamma(X,  \mathcal{O} ((\ell+r)D_0))$
is an irreducible highest weight module with highest weight
\[
\lambda_r(\ell)=\begin{cases} s_r\cdots
s_2s_1(\ell\varpi_1+\rho)-\rho,
\text{ if } \ell+r\geq 0 \text{ and } r\leq n,\\
s_{r-1}\cdots s_2s_1(\ell\varpi_1+\rho)-\rho, \text{ if } \ell+r< 0
\text{ and } r\geq 1,\\
0,\text{ if } \ell+r=0 \text{ and } r=n+1.
\end{cases}
\]
In the remaining cases (namely $\ell+r\geq 0, r=n+1$ and $\ell+r<0,
r=0$), $\Gamma(X,\mathcal{O}(\ell+r)D_0)=0$.
\end{Thm}

\begin{proof}
We use the description \eqref{e-Richard Parker} of the sheaf of
sections of a line bundle with divisor $D=\sum a_i D_i$. The exact
sequence \eqref{e-xs2} is
\[
0\to K\to \mathbb Z^{n+1}\stackrel p\to\sum_{i=1}^n \mathbb Z k_i\to
0,
\]
The kernel $K$ of the map $p\colon \epsilon_i\mapsto k_i$ is spanned
by $\sum_{i\notin I}\epsilon_i-\sum_{i\in I}\epsilon_i$. The
character $\chi$ is (the restriction to $T_K$ of) the character
$z\mapsto z^a=z_1^{a_1}\cdots z_{n+1}^{a_{n+1}}$ of
$(\K^\times)^{n+1}$. The sections on the open set $U_0$ (the open
orbit of $T_N$) are then spanned by the monomials $Q^{\mu+a}$ where
$Q_i$ are the coordinates on $Y\subset\K^{n+1}$ and $\mu$ runs over
$K^\perp=\{\mu\in\mathbb Z^{n+1}\,|\,\sum_{i\notin
I}\mu_i=\sum_{i\in I}\mu_i\}$. The sections $Q^{\mu+a}$ extend to
global regular sections if and only if $\mu_i+a_i\geq 0$ for all
$i$. Let $\sigma_I\colon \mathbb Z^{n+1}\to\mathbb Z^{n+1}$ be the
involution which changes the sign of the coordinates labeled by $I$.
Then we get the following (well-known) description of $\Gamma(X,
\mathcal{O} (\sum a_iD_i))$: a basis is
\begin{equation}\label{e-Crouch}
\{Q^{\nu}\,|\,\sum_{i\notin I}\nu_i-\sum_{i\in I}\nu_i=m,\,
\nu_i\geq0,\, i=1,\dots,n+1\},
\end{equation}
where $ m=\sum_{i\notin I}a_i-\sum_{i\in I}a_i$. This basis consists
of weight vectors. Let us compute their weight in the weight lattice
$P=\mathbb Z^{n+1}/\mathbb Z\cdot (1,\dots,1)$:
$x=\mathrm{diag}(x_1,\dots,x_{n+1})\in\mathfrak h$ acts as $F_I(\sum
x_iQ_iP_i)=\sum_{i\notin I}x_iQ_iP_i-\sum_{i\in I}x_iP_iQ_i$, thus
\[
x\cdot Q^{\nu}=\left(\sum_{i\notin I}\nu_i x_i-\sum_{i\in
I}(\nu_i+1)x_i\right)Q^{\nu}.
\]
Thus $Q^{\nu}$ has weight $\lambda=\sigma_I(\nu)-\sum_{i\in
I}\epsilon_i$ modulo $\mathbb Z\cdot(1,\dots,1)$. Now suppose that
the divisor $D=\phi_I(D^0)$ is obtained from a divisor $D^0=\sum
a^0_iD_i$ on $\mathbb P^n$, so that $a=\sigma_I(a^0)-\sum_{i\in
I}\epsilon_i$. Let $\ell=\sum a_i^0$ so that $m=\ell+r$. Then the
weights appearing in $\Gamma(X,  \mathcal{O} ( D))$ are the classes
in $P$ of the integer vectors $\lambda$ such that
\[
\sum_{i=1}^{n+1}\lambda_i=\ell,\qquad \lambda_i\geq 0 \text{ if }
i\notin I,\quad\lambda_i\leq -1 \text{ if } i\in I,
\]
where $\ell=\sum_ia_i^0=m-|I|$. The corresponding weight spaces are
one-dimensional, spanned by $Q^{\nu}$ with
$\nu=\sigma_I(\lambda)-\sum_{i\in I}\epsilon_i$.

Let now $I=\{1,\dots,r\}$. We need to identify the primitive weight
vectors in $\Gamma(X,  \mathcal{O} (D))$, namely the monomials $Q^\nu$ annihilated
by the action of $\mathfrak n_+$.

\begin{Lem}\label{l-Tomatlan}
The only primitive monomial in $\Gamma(X,  \mathcal{O} (D))$ is $1$
if $m=0$, $Q_r^{-m}$ if $m<0$, $r\geq1$ and $Q_{r+1}^{m}$ if $m>0,
r\leq n$. If $m<0$ and $r=0$ or if $m>0$ and $r=n+1$ then
$\Gamma(X,\mathcal{O}(D))=0$.
\end{Lem}

\begin{proof}
The monomial $Q^\nu$ is primitive if and only if $e_i\cdot Q^\nu=0$
for the Chevalley generators $e_i$, $i=1,\dots, n$, associated with
simple roots. These generators act as
\[
e_i=F_I(Q_iP_{i+1})=
 \begin{cases}-P_iQ_{i+1}, & i<r,\\
 P_iP_{i+1},&i=r,\\
 Q_iP_{i+1},&i>r,\\
 \end{cases}
\]
The only solutions of $e_i\cdot Q^\nu=0$ are powers of $Q_r$ and of
$Q_{r+1}$. The condition on the sum of exponents in \eqref{e-Crouch}
implies that the exponents are $-m$ and $m$, respectively. The last
assertion follows from the fact that there is no homogeneous
polynomial of negative degree.
\end{proof}
If $m=\ell+r\geq0$ then the weight of $Q_{r+1}^m$ is
\[
\lambda_r(\ell)=\sigma_I((\ell+r)\epsilon_{r+1})-\sum_{i\in
I}\epsilon_i=(\ell+r)\epsilon_{r+1}-\sum_{i=1}^r\epsilon_i.
\]
If $m<0$ then the weight of $Q_r^m$ is
\[
\lambda_r(\ell)=\sigma_I(-(\ell+r)\epsilon_r)-\sum_{i\in
I}\epsilon_i=(\ell+r-1)\epsilon_r-\sum_{i=1}^{r-1}\epsilon_i.
\]
To identify these weights with the vector in the $\rho$-shifted Weyl
orbit of $\ell \varpi_1$, recall that $\varpi_1=\epsilon_1$, that
$s_i$ is transposition of $\epsilon_i$ and $\epsilon_{i+1}$ and that
$s_i\rho-\rho=-\alpha_i=\epsilon_{i+1}-\epsilon_{i}$. It follows
that there is a nonzero module homomorphism $M(\lambda_r(\ell))\to
\Gamma(X,  \mathcal{O} (D))$ from the Verma module with highest
weight $\lambda_r(\ell)$. The image is irreducible by Lemma
\ref{l-Tomatlan}, since any proper submodule would contain a
primitive monomial distinct from the highest weight vectors. It
remains to show that this homomorphism is surjective. This follows
from the next Lemma.

\begin{Lem}
The $\mathcal D_{  \mathcal{O}(D)}(X)$-module $\Gamma(X, \mathcal{O}
(D))$ (when $\neq0$) is generated by the primitive monomial of Lemma
\ref{l-Tomatlan}.
\end{Lem}
\begin{proof} By Musson's theorem, $\mathcal D_{\mathcal{O}(D)}(X)$ is spanned by the images of the
monomial differential operators $Q^\lambda P^\mu$ such that $\tau
=\lambda-\mu$ obeys $\sum_{i\notin I}\tau_i=\sum_{i\in I}\tau_i$. If
$m\geq 0$, a general monomial $Q^\nu$ of the basis \eqref{e-Crouch}
can be obtained from the highest weight vector by the differential
operator $Q^\nu P_{r+1}^m/m!$ on the highest weight vector
$Q_{r+1}^m$. If $m<0$, take $Q^\nu P_{r}^m/m!$.
\end{proof}
This concludes the proof of Theorem \ref{t-5}.
\end{proof}

\begin{Rem}
The symmetric group $S_{n+1}$ of $\mathfrak{sl}_{n+1}$ acts on
$U(\mathfrak{sl}_{n+1})$ by automorphisms (it is the Weyl group) and
on $\mathbb P^n$ and $\mathcal O(\ell)$ by permutations of
homogeneous coordinates. The map
$U(\mathfrak{sl}_{n+1})\to\mathcal{D}_{\mathcal O(\ell)}(\mathbb P^n)$ is
$S_{n+1}$-equivariant. The modules $\Gamma(X,\mathcal
O((\ell+r)D_0))$ corresponding to the other subsets $I$ are related
by the Weyl automorphism of $U(\mathfrak{sl}_{n+1})$ associated with
any permutation sending $I$ to $\{1,\dots,r\}$ to the modules of the
Theorem. They are thus highest weight modules for other Cartan
decompositions.
\end{Rem}
\begin{Rem}
In the case of the variety of rank one matrices $Z_r$, $2\leq r\leq
n-1$, only the trivial vector bundle extends to the closure $\bar
Z_r$ in the space of all matrices, considered in
\cite{Levasseuretal}. In this case we recover the result of
\cite{Levasseuretal} that $\Gamma(Z_r,\mathcal O)\simeq
L(-\omega_r)$.
\end{Rem}

\subsection{Higher cohomology groups}
Finally, we want to discuss the $\mathfrak{sl}_{n+1}$-module structure of higher cohomology groups for $X$ a regular toric variety of the kind considered in (c) and (e) in Subsection~\ref{ss-31}, related to $\mathbb P^n$ by an $I$-reflection, with $I=\{1,\dots,r\}$ and $2\leq 2\leq n$: in particular, we consider $X=\tilde Z_r^+$, for $2\leq r\leq n-1$, or $X=\widetilde{\mathbb A^n}$, for $r=n$.
Let as above $D_0$ be a generator of the Picard group of $X$; for the subsequent computations, assume $D_0=D_{n+1}$.
Finally, let $m=$ be an integer number.
\begin{Thm}\label{thm-coh}
For $I$, $X$, $D_0$ and $m$ as above, we have
\[
H^*(X,\mathcal O(m D_0))\cong \bigoplus_{\lambda\in \mathbb N^{n-r+1}} H^*(\mathbb P^{r-1},\mathcal O(-m+|\lambda|)).
\]
\end{Thm}
\begin{proof}
We compute the cohomology of $X$ with values in $\mathcal O(m D_0)$ by means of \v Cech cohomology.
For this purpose, recall that $X$ is a torus quotient and that it has an open covering by sets $U_i=\{Q_i\neq 0\}$, $i=1,\dots,r$, and $Q_i$ denoting affine coordinates on $\mathbb A^{n+1}$.
Further, the Weil divisor $m D_0$ has a Cartier representative of the form $f_{i,m}(Q)=Q_i^m Q_{n+1}^m$, $i=1,\dots,r$: the transition functions of the corresponding line bundle $\mathcal O(m D_0)$ have the form $(Q_i/Q_j)^m$ on $U_i\cap U_j$, $1\leq i\neq j\leq r$.
Finally, the variety  $X$ admits local affine coordinates over $U_i$, namely 
\[
[(Q_1,\dots,Q_r,Q_{r+1},\dots,Q_{n+1})]\mapsto \left(\frac{Q_1}{Q_i},\dots,\frac{Q_r}{Q_i},Q_iQ_{r+1},\dots,Q_iQ_{n+1}\right),
\]
where we omit the $i$-th entry.
We denote by $(z_1,\dots,z_{r-1},w_1,\dots,w_{n-r+1})$ the corresponding local coordinates on $\mathbb A^n$.

W.r.t.\ the open covering $\{U_i\}$, a \v Cech $p$-cochain on $X$ with values in $\mathcal O(m D_0)$ is represented by a family $\sigma_{(i_0,\dots,i_p)}$ of regular sections of $\mathcal O(m D_0)$ on the intersection $U_{i_0}\cap\cdots \cap U_{i_p}$.
E.g.\ w.r.t.\ the affine coordinates over $U_{i_0}$, the component $\sigma_{(i_0,\dots,i_p)}$ is represented by a regular function $P_{(i_0,\dots,i_p)}(z,w)$ on $\bigcap_{i=1}^p \{z_i\neq 0\}$, admitting an expansion of the form
\[
P_{(i_0,\dots,i_p)}(z,w)=\sum_{\lambda \in\mathbb N^{n-r+1}}P_{(i_0,\dots,i_p)}^\lambda(z)w^\lambda,
\] 
where $P_{(i_0,\dots,i_p)}^\lambda$ are now regular functions on $\bigcap_{i=1}^p\{z_i\neq 0\}$.

The cocycle- and coboundary-conditions can now be rewritten in terms of the previous expansion w.r.t.\ $w$, using the previous transition functions for $\mathcal O(m D_0)$: denoting the degree of $\lambda\in \mathbb N^{n-r+1}$ by $|\lambda|$, since neither a change of coordinates on $X$ nor a change of trivialization on $\mathcal O(m D_0)$ affects the degree $|\lambda|$ of $w$, the cocycle- and coboundary conditions for \v Cech cohomology on $X$ can be rewritten as cocycle- and coboundary conditions for $P_{(i_0,\dots,i_p)}^\lambda$, viewed as sections on a Serre bundle $\mathcal O(-m+|\lambda|)$ over $\mathbb P^{r-1}$, for any multiindex $\lambda$ as above.
Notice that the transition functions of $\mathcal O(m D_0)$ depend only on the homogeneous coordinates $Q_i$, for $i=1,\dots,r$, and when $Q_i=0$, $i=r+1,\dots,n+1$, they correspond to the transition functions of the Serre bundle $\mathcal O(-m)$.
Finally, notice that the affine coordinates $w_i$ are all multiplied by some coordinate $z_j$, when performing a corresponding coordinate change on $\mathbb P^{r-1}$, whence the shift by the degree $|\lambda|$ in the Serre bundles.

Summarizing, any \v Cech cocycle on $X$ of degree $p$ with values in $\mathcal O(m D_0)$ is equivalent to a family of \v Cech cocycles of the same degree on $\mathbb P^{r-1}$ with values in certain Serre bundles related to $\mathcal O(-m)$ by shifts of the parameter; the same holds true for \v Cech coboundaries, and these two facts yield the above isomorphism.
\end{proof}
In particular, since the cohomology of projective spaces $\mathbb P^{r-1}$ is concentrated in degree $0$ and $r-1$, it follows also that the cohomology of $X$ with values in $\mathcal O(m D_0)$ is concentrated in degree $0$ and $r-1$.

We can also describe explicitly the isomorphism of Theorem~\ref{thm-coh} using homogeneous coordinates $Q_i$ on $X$, knowing the description of the cohomology of projective spaces.
Namely, if $k\geq 0$, the $0$-th cohomology of $\mathbb P^{r-1}$ with values in $\mathcal O(k)$ is generated by monomials $Q_1^{\mu_1}\cdots Q_r^{\mu_r}$, $\mu_i\geq 0$ and $\sum_{i=1}^r \mu_i=k$, while, if $k\leq -r$, the $r-1$-th cohomology with values in $\mathcal O(k)$ is generated (modulo coboundaries) by monomials of the same kind, with $\mu_i<0$ and $\sum_{i=1}^r \mu_i=k$ (in all other cases, the cohomologies are trivial).

Then, the isomorphism of Theorem~\ref{thm-coh} can be written as
\begin{multline*}
H^*(\mathbb P^{r-1},\mathcal O(-m+|\lambda|))\ni Q_1^{\mu_1}\cdots Q_r^{\mu_r}\mapsto \\\mapsto Q_1^{\mu_1}\cdots Q_r^{\mu_r}Q_{r+1}^{\lambda_{r+1}}\cdots Q_{n+1}^{\lambda_{n+1}}\in H^*(X,\mathcal O(m D_0)),
\end{multline*}
for any $\lambda\in \mathbb N^{n-r+1}$.
In particular, a basis of the cohomology of $X$ with values in $\mathcal O((l+r)D_0)$ is given $i)$ (in degree $0$) by monomials $Q^\nu$, with $\sum_{i=1}^r \nu_i+m=\sum_{i=r+1}^{n+1}\nu_i$ and $\nu_i\geq 0$, $i=1,\dots,n+1$, (see also Subsection 4.2.) and $ii)$ (in degree $r-1$) by monomials $Q^\nu$, with $\sum_{i=1}^r \nu_i+m=\sum_{i=r+1}^{n+1}\nu_i$ and $\nu_i<0$, for $i=1,\dots,r$, and $\nu_i\geq 0$, for $i=r+1,\dots,n+1$.
It follows immediately that the $0$-th cohomology of $X$ with values in $\mathcal O(m D_0)$ is always infinite-dimensional and non-trivial, while the $r-1$-th cohomology is always finite-dimensional and is non-trivial exactly when $m\geq r$.

The $\mathfrak{sl}_{n+1}$-module structure on the $0$-th cohomology was discussed in Subsection 4.2.: we now consider the cohomology of degree $r-1$.
We set now $m=\ell+r$ for an integer $\ell$; then, by Theorem~\ref{thm-coh} and the discussion thereafter, we need only discuss the case $l\geq 0$.
\begin{Thm}\label{thm-irr}
For $X$, $I$, $D_0$ as in Theorem~\ref{thm-coh}, and $\ell\geq 0$, there is an isomorphism 
\[
H^{r-1}(X,\mathcal O((\ell+r) D_0))\cong L(\ell \varpi_1)
\]
of $\mathfrak{sl}_{n+1}$-modules.
\end{Thm}  
\begin{proof}
We first compute the possible primitive vectors in $H^{r-1}(X,\mathcal O((\ell+r) D_0))$; for this purpose, recall the expressions for the Chevalley generators $e_i$ in Lemma~\ref{l-Tomatlan}.
A monomial $Q^\nu$ representing a non-trivial class in $H^{r-1}(X,\mathcal O((\ell+r) D_0))$ is annihilated by $e_i$, $i=1,\dots,n$, if and only if $\nu_i=0$, for $i=r+1,\dots,n+1$, by direct computations, using the generators $e_i$, $i=r,\dots,n+1$.
Thus, $Q^\nu$ has the form $Q_1^{\nu_1}\cdots Q_r^{\nu_r}$, with $\nu_i<0$, and $\sum_{i=1}^r\nu_i=-\ell-r$.
Using the remaining generators $e_i$, we see that $e_i Q^\mu$ vanishes in cohomology if and only if $\nu_i=-1$, $i=2,\dots,r$; then, the remaining exponent is automatically $\nu_1=-\ell-1$.
Thus, for any $\ell\geq 0$, the only primitive vector in $H^{r-1}(X,\mathcal O((\ell+r) D_0))$ is the monomial $Q_1^{-\ell-1}Q_2^{-1}\cdots Q_r^{-1}$.

The primitive vector $Q_1^{-\ell-1}Q_2^{-1}\cdots Q_r^{-1}$ is also a weight vector: its weight is readily computed, since
\[
x\cdot(Q_1^{-\ell-1}Q_2^{-1}\cdots Q_r^{-1})=
\ell x_1 (Q_1^{-\ell-1}Q_2^{-1}\cdots Q_r^{-1}),
\]
for any $x=\mathrm{diag}(x_1,\dots,x_{n+1})\in \mathfrak{sl}_{n+1}$, which acts in this situation as the differential operator $-\sum_{i=1}^r x_i Q_iP_i+\sum_{i=r+1}^{n+1} x_i Q_iP_i-\sum_{i=1}^r x_i$ (see also the proof of Theorem~\ref{t-5}).
Therefore, $H^{r-1}(X,\mathcal O((\ell+r) D_0))$ has a unique primitive vector of weight $\ell \varpi_1$.

Moreover, the module $H^{r-1}(X,\mathcal O((\ell+r) D_0))$ is generated by the primitive monomial $Q_1^{-\ell-1}Q_2^{-1}\cdots Q_r^{-1}$.
Namely, a general differential operator on $\mathcal D_{\mathcal O((\ell+r) D_0)}(X)$ is of the form $Q^\lambda P^\mu$, with $\lambda$, $\mu$ positive multiindices, such that $\tau=\lambda-\mu$ satisfies $\sum_{i=1}^r\tau_i=\sum_{i=r+1}^{n+1} \tau_i$.
On the other hand, $Q^\nu$ is a generator of $H^{r-1}(X,\mathcal O((\ell+r) D_0))$ if $\nu_i<0$, $i=1,\dots,r$, $\nu_i\geq 0$, $i=r+1,\dots,n+1$, and $\sum_{i=1}^r \nu_i+\ell+r=\sum_{i=r+1}^{n+1}\nu_i$.
For such a weight $\nu$, the differential operator 
\[
D_{\nu}=\frac{(-1)^{\sum_{i=1}^r \nu_i+r}}{\prod_{i=1}^r (-(\nu_i+1))!}Q_1^\ell Q_{r+1}^{\nu_{r+1}}\cdots Q_{n+1}^{\nu_{n+1}}P_1^{-(\nu_1+1)}\cdots P_r^{-(\nu_r+1)} 
\]
belongs, by a direct computation, to $\mathcal D_{\mathcal O((\ell+r) D_0)}(X)$, and it is easy to verify that 
\[
Q^\nu=D_\nu(Q_1^{-\ell-1}Q_2^{-1}\cdots Q_r^{-1}).
\]
All these computations show, by the very definition of Verma modules, that there is a surjective, non-trivial module homomorphism $M(l\varpi_1)\to H^{r-1}(X,\mathcal O((\ell+r) D_0))$, whose image is irreducible (since there is only one primitive vector up to multiplication by $\mathbb C^\times$): this yields the above isomorphism.  
\end{proof}
Notice that the highest weight $l\varpi_1$ is dominant, which re-proves the fact that $H^{r-1}(X,\mathcal O((\ell+r) D_0))$, when non-trivial, is always finite-dimensional.

\begin{Rem}\label{rem-chev}
Theorem~\ref{thm-irr} refers only to $X=\tilde Z_r^{+}$ or $X=\widetilde{\mathbb A^n}$, with $2\leq r\leq n$: the reason is that $i)$ the resolution of singularities $\tilde Z_r^-$, corresponding to $2\leq r\leq n-1$, are related to $\tilde Z_{n-r}^+$ by an isomorphism of the corresponding fans, which lifts to the Chevalley involution of $\mathfrak{sl}_{n+1}$ on the corresponding (twisted) rings of differential operators, and $ii)$, if $r=1$, $X=\widetilde{\mathbb A^n}$, and the two blow-ups $r=1$ and $r=n$ are related to each other by an automorphism of the corresponding fan, which also lifts to the Chevalley involution at the level of differential operators.
Hence, by means of the Chevalley involution, we can deduce the $\mathfrak{sl}_{n+1}$-module structure of the twisted cohomologies of $\tilde Z_r^-$, $2\leq r\leq n-1$, and of $\widetilde{\mathbb A^n}$, for $r=1$, from Theorem~\ref{thm-coh} and~\ref{thm-irr}.
\end{Rem}


\begin{bibdiv}
\begin{biblist}
\bib{Audin}{book}{
   author={Audin, Mich{\`e}le},
   title={The topology of torus actions on symplectic manifolds},
   series={Progress in Mathematics},
   volume={93},
   note={Translated from the French by the author},
   publisher={Birkh\"auser Verlag},
   place={Basel},
   date={1991},
   pages={181},
   isbn={3-7643-2602-6},
   review={\MR{1106194 (92m:57046)}},
}

\bib{BerestWilson}{article}{
  author={Berest, Yuri},
  author={Wilson, George},
  title={Differential Isomorphism and Equivalence of Algebraic Varieties},
  eprint={http://arxiv.org/abs/math/0304320v1},
  date={2003}
}


\bib{BorhoBrylinski}{article}{
   author={Borho, Walter},
   author={Brylinski, Jean-Luc},
   title={Differential operators on homogeneous spaces. I. Irreducibility of
   the associated variety for annihilators of induced modules},
   journal={Invent. Math.},
   volume={69},
   date={1982},
   number={3},
   pages={437--476},
   issn={0020-9910},
   review={\MR{679767 (84b:17007)}},
}

\bib{Cox}{article}{
   author={Cox, David A.},
   title={The homogeneous coordinate ring of a toric variety},
   journal={J. Algebraic Geom.},
   volume={4},
   date={1995},
   number={1},
   pages={17--50},
   issn={1056-3911},
   review={\MR{1299003 (95i:14046)}},
}


\bib{Fulton}{book}{
   author={Fulton, William},
   title={Introduction to toric varieties},
   series={Annals of Mathematics Studies},
   volume={131},
   note={;
   The William H. Roever Lectures in Geometry},
   publisher={Princeton University Press},
   place={Princeton, NJ},
   date={1993},
   pages={xii+157},
   isbn={0-691-00049-2},
   review={\MR{1234037 (94g:14028)}},
}
\bib{Goncharov}{article}{
   author={Goncharov, A. B.},
   title={Constructions of Weyl representations of some simple Lie algebras},
   language={Russian},
   journal={Funktsional. Anal. i Prilozhen.},
   volume={16},
   date={1982},
   number={2},
   pages={70--71},
   issn={0374-1990},
   review={\MR{659169 (84d:22023)}},
}
\bib{Jantzen}{book}{
   author={Jantzen, Jens Carsten},
   title={Einh\"ullende Algebren halbeinfacher Lie-Algebren},
   language={German},
   series={Ergebnisse der Mathematik und ihrer Grenzgebiete (3) [Results in
   Mathematics and Related Areas (3)]},
   volume={3},
   publisher={Springer-Verlag},
   place={Berlin},
   date={1983},
   pages={ii+298},
   isbn={3-540-12178-1},
   review={\MR{721170 (86c:17011)}},
}

\bib{Levasseuretal}{article}{
   author={Levasseur, T.},
   author={Smith, S. P.},
   author={Stafford, J. T.},
   title={The minimal nilpotent orbit, the Joseph ideal, and differential
   operators},
   journal={J. Algebra},
   volume={116},
   date={1988},
   number={2},
   pages={480--501},
   issn={0021-8693},
   review={\MR{953165 (89k:17028)}},
}
\bib{Musson}{article}{
   author={Musson, Ian M.},
   title={Differential operators on toric varieties},
   journal={J. Pure Appl. Algebra},
   volume={95},
   date={1994},
   number={3},
   pages={303--315},
   issn={0022-4049},
   review={\MR{1295963 (95i:16026)}},
}
\bib{Musson2}{article}{
   author={Musson, Ian M.},
   title={Actions of tori on Weyl algebras},
   journal={Comm. Algebra},
   volume={16},
   date={1988},
   number={1},
   pages={139--148},
   issn={0092-7872},
   review={\MR{921946 (88k:17013)}},
}
\bib{Musson-Rueda}{article}{
    author={Musson, Ian M.},
    author={Rueda, Sonia L.},
    title={Finite dimensional representations of invariant differential
    operators},
    journal={Trans. Amer. Math. Soc.},
    volume={357},
    date={2005},
    number={7},
    pages={2739--2752 (electronic)},
    issn={0002-9947},
    review={\MR{2139525 (2006a:16034)}},
 }

                
\end{biblist}
\end{bibdiv}
\end{document}